\documentclass[11pt]{article}

\usepackage{amsfonts}
\usepackage{amsmath}
\usepackage{amssymb}
\usepackage{amsthm}

\def\calP{\mathcal{P}}

\def\calpd{\calP_{p,\delta}}
\def\calpp{\calP_p}

\def\hcalP{\hat{\mathcal{P}}}
\def\hcaltat{\hcalP_{\tau,t}}
\def\hcaltaiti{\hcalP_{\tau_i,t_i}}
\def\hN{\hat{N}}

\def\tL{\tilde{L}}
\def\ra{\rightarrow}
\def\lra{\leftrightarrow}
\def\da{\downarrow}
\def \Z {\mathbb Z}
\def \R {\mathbb R}
\def \N {\mathbb N}

\def\al{\alpha}
\def\th{\theta}
\def\hth{\hat{\theta}}
\def\heta{\hat{\eta}}
\def\be{\beta}
\def\ep{\varepsilon}
\def\de{\delta}

\def\si{\sigma}
\def\de{{\delta}}

\def\om{\omega}
\def\Om{\Omega}

\newtheorem{thm}{Theorem}[section]
\newtheorem{prop}[thm]{Proposition}
\newtheorem{lem}[thm]{Lemma}
\newtheorem{cor}[thm]{Corollary}

\theoremstyle{plain}

\newtheorem{rem}[thm]{Remark}

\begin{document}

\title{Continuity for self-destructive percolation in the plane}

\author{J. van den Berg\footnote{Part of vdB's research has been funded by the Dutch BSIK/BRICKS project.} \,\, R. Brouwer  and
B. V\'{a}gv\"{o}lgyi  \\
{\small CWI and VUA} \\
{\footnotesize email: J.van.den.Berg@cwi.nl; Rachel.Brouwer@cwi.nl; B.Vagvolgyi@few.vu.nl}
}
\date{}
\maketitle

\begin{abstract}
A few years ago (see \cite{BeB}) two of us introduced, motivated by the study of certain forest-fire
processes, the self-destructive percolation model (abbreviated as sdp model).
A typical configuration for the
sdp model with parameters $p$ and $\de$ is generated in
three steps: First we generate a typical configuration for the ordinary percolation model
with parameter $p$. Next, we make all sites in the infinite occupied cluster vacant. 
Finally, each site
that was already vacant in the beginning or made vacant by the above action,
becomes occupied with probability $\de$ (independent of the other sites).

Let $\th(p,\de)$ be the probability that some specified vertex belongs,
in the final configuration, to an infinite occupied cluster. 
In our earlier paper we stated the conjecture that, for the square lattice and other planar lattices,
the function $\th(\cdot,\cdot)$ has a discontinuity at points of the form $(p_c, \delta)$, with $\de$
sufficiently small. We also showed (see \cite{BeB2}) remarkable
consequences for the forest-fire models.

The conjecture naturally raises the question whether the function $\th(\cdot,\cdot)$ is continuous outside some region
of the above mentioned form.
We prove that this is indeed the case. An important ingredient in our proof is 
a (somewhat stronger form of a) recent ingenious RSW-like percolation result of
Bollob\'{a}s and Riordan (\cite{BoR}).

\end{abstract}

\begin{section}{Introduction and outline of results}
\begin{subsection}{Background and motivation}
The self-destructive percolation model on the square lattice is described as follows:
First we perform independent
site percolation on this lattice: 
we declare each site \emph{occupied} with probability $p$, and \emph{vacant}
with probability $1-p$, independent of the other sites. We will use the notation $\{V \lra W\}$ for the event that
there is an occupied path from the set of
sites $V$ to the set of sites $W$.
We write $\{V \lra \infty\}$ for the event that there is an infinite occupied path starting at $V$. \\
Let, as usual, $\th(p)$ denote the probability that a given site, say $O = (0,0)$,
belongs to an infinite occupied cluster.
It is known that there is a critical value $0 < p_c <1$ such that $\th(p) > 0$ for all
$p > p_c$, and $\th(p) = 0$ for all $p \leq p_c$.
%(This is also known for several other lattices including the triangular and the honeycomb
%lattice and regular trees; it is strongly believed that this also holds for the $d$-dimensional
%hypercubic lattice, $d \geq 3$, but this has only been proved for sufficiently high $d$).
Now suppose that, by some catastrophe, the infinite occupied cluster (if present) is destroyed;
that is, each site in this cluster becomes vacant.
Further suppose that after this catastrophe we give the sites independent
`enhancements', as follows: Each site that was already vacant in the beginning, or was {\it made} vacant by the
catastrophe,
becomes occupied with probability $\de$, independent of the others.
%Those sites that were already occupied after
%the catastrophe remain occupied.
Let $\calpd$ be the distribution of the final configuration.

A more formal, and often very convenient description
of the model is as follows:
Let $X_i$, $i \in \Z^2$ be independent $0-1$ valued random variables,
each $X_i$ being $1$ with probability $p$ and $0$ with probability $1-p$.
Further,
let $Y_i$, $i \in \Z^2$, be independent $0-1$ valued random variables,
each $Y_i$ being $1$ with probability $\de$ and $0$ with probability $1-\de$.
Moreover, we take the collection of $Y_i$'s independent of that of the $X_i$'s.
Let $X_i^* , i \in \Z^2$ be defined by

\begin{eqnarray}
X_i^* &=&
\left\{ \begin{array}{ll}
& 1 \mbox{ if } X_i = 1 \mbox { and there is no $X$- occupied path from } i \mbox{ to }
\infty \\
& 0 \mbox{ otherwise},
\end{array}
\right.\nonumber \\
\end{eqnarray}
where by `$X$-occupied path' we mean a path on which each site $j$ has $X_j = 1$.
Finally, define $Z_i = X_i^* \vee Y_i$.
This collection $(Z_i, i \in \Z^2)$ is (with 0 meaning `vacant' and 1 `occupied') what we called `the final
configuration', and 
the above mentioned $\calpd$ is its distribution.

We use the notation $\th(p, \de)$ for the probability that, in the final
configuration, $O$ is in an infinite occupied cluster:

$$\th(p,\de) := \calpd(O \lra \infty).$$

Note that
$O$ is occupied in the final configuration if and only if the above
mentioned enhancement was successful, or $O$ belonged initially
(before the catastrophe) to a non-empty but finite occupied cluster.
This gives
$$\calpd(O \mbox{ is occupied }) = \de + (1-\de) (p -\th(p)).$$

Also note that, in the case that $p \leq p_c$, nothing happens in the above
catastrophe, so that in the final configuration the
sites are independently occupied with probability $p +(1-p) \de$.
Formally, if $p \leq p_c$, then

\begin{equation}
\calpd = \calP_{p+(1-p)\de},
\label{ordinary}
\end{equation}

\noindent
where we use the notation $\calP_p$ for the product measure with
parameter $p$.
In particular,

\begin{equation}
\label{thpcd}
\th(p_c, \de) = \th(p_c + (1-p_c)\de) > 0,
\end{equation}
for each $\de > 0$.

\begin{rem}
Most of what we said above has straightforward analogs for arbitrary countable graphs, but there are subtle 
differences.  For instance, on the cubic lattice it has not yet been proved that $\th(p_c) = 0$ (although this 
is generally believed to be true). So, for that lattice, \eqref{ordinary} with $p = p_c$, and hence \eqref{thpcd},
are not rigorously known.
\end{rem}

\smallskip\noindent
It is also clear from the construction that $\calpd$ stochastically dominates
$\calP_{\delta}$. Hence, if $\de > p_c$ then
$\th(p,\de) \geq \th(\de) > 0$ for all $p$.

\medskip\noindent
It turns out (see Proposition 3.1 of \cite{BeB}) 
that, if $p>p_c$, a `non-negligible'
enhancement is needed after the catastrophe to create again an infinite occupied cluster.
More precisely, for each $p > p_c$ there is a $\de > 0$ with $\th(p, \de) = 0$.
A much more difficult question is whether the needed enhancement
goes to $0$ as $p \da p_c$.
By \eqref{thpcd} one might be tempted to reason intuitively that this is indeed the case.
In \cite{BeB} it was shown that for the analogous model on the binary tree this is correct.
However, in \cite{BeB} a conjecture is presented which says, in particular, that for the square lattice (and other planar lattices)
there is a $\de > 0$ for which $\th(p,\de) = 0$ for all $p > p_c$. 
%For that model we prove in Section 5 that for
%all $\delta >0,$
%$$\lim_{p \da p_c} \th(p, \de) = \th(p_c, \de) >0.$$ \\
In Section 4 of \cite{BeB} and in \cite{BeB2} we showed remarkable consequences for certain forest-fire models.

\smallskip\noindent
Note that, since $\th(p_c,\de) > 0$, the above conjecture says that the function $\th(\cdot,\cdot)$ has
discontinuities at points of the form $(p_c,\de)$ with $\de$ sufficiently small.
This naturally raises the question whether this function is continuous in the complement of a region of such form:
is there a
$\de > 0$ such that $\th(\cdot,\cdot)$ is continuous outside the set $\{p_c\} \times [0,\de]$?
In the next subsection we state that this is indeed the case, and give a summary of the methods and intermediate results used
in the proof. At the end of Section \ref{sec6}  we point out why our proof does not work at points $(p_c,\de)$ with
small $\de$. We hope our arguments provide a better understanding of the earlier mentioned conjecture and will trigger
new attempts to prove (or disprove) it.
\end{subsection}

\begin{subsection}{Outline of results}
The conjecture mentioned in the previous subsection 
raises the natural question whether $\th(\cdot,\cdot)$ {\it is} continuous outside
the indicated `suspected' region. The following theorem states that this is indeed the case.
\begin{thm}
There is a $\de \in (0,1)$ such that the function $\th(\cdot,\cdot)$ is continuous outside
the segment $\{p_c\} \times (0, \de)$.
\label{mainthm}
\end{thm}

As could be expected, the proof widely uses tools and results from ordinary percolation. However, the 
dependencies introduced by the self-destructive mechanism cause complications. Until
recently, a serious obstacle was the absence of a suitable RSW-like theorem.
This obstacle could be removed by the use of (a modified and somewhat stronger form of) a recent theorem
of Bollob\'{a}s and Riordan (\cite{BoR}).

A rough outline of the proof of Theorem \ref{mainthm}, and the needed intermediate results that are interesting in
themselves, is as follows:
In section 2 we list some basic properties of our model, which will be used later.
The results in Section 3, which are also contained in the recent PhD thesis \cite{Bro} of one of us, show
that if $\th(\cdot, \cdot)$ is strictly
positive in some open region, then it is continuous on this region.
It is also shown that if $\th(p,\de) = 0$, then $\th(\cdot,\cdot)$ is continuous at $(p,\de)$.
These two results reduce the proof of Theorem \ref{mainthm} to showing that if $\th(p,\de) > 0$ and $p \neq p_c$, then
$\th(p,\de) > 0$ in an open neighborhood of $(p, \de)$. 
This in turn requires a suitable finite-size criterion (see below) for sdp.
In Section 4 we give the modified form of the Bollob\'{a}s-Riordan theorem.
%The result by Bollob\'{a}s and Riordan was intended for the Voronoi percolation model. We first give a slightly stronger
%version of their result. We also observed that, with minor changes, their proof-by-contradiction could be turned
%into a (in our view) somehat more `direct' proof. This more direct setup gives (in our view) some more clarity and 
%leads immediately to the slightly stronger form we needed. We follow closely the list of claims in
%the Bollob\'{a}s-Riordan proof
%and indicate where we make changes. This is still for the Voronoi model. As Bollob\'{a}s and Riordan argued at the end of their
%proof, their result holds not only for Vornoi percolation but for a wide class of models with certain properties. This 
%is also the case for our slight modified form and we show that the sdp model has the  required properties.
This is used in Section 5 to obtain the above mentioned finite-size criterion.
Finally, in Section 6 we combine these results and prove the main theorem.

We end this section with the following remark: When we say that a function $f$ is `increasing' (`decreasing') this should,
unless this is preceded by the word `strictly',
be interpreted in the weak sense: $x < y$ implies $f(x) \leq  f(y)$.

\end{subsection}
\end{section}

\begin{section}{Basic properties}
In this section we state some basic properties which will be used later.

First some more terminology and notation: 
If $v = (v_1,v_2)$ and $w= (w_1,w_2)$ are two vertices, we let
$|v-w|$  denote their
(graph) distance $|v_1 - w_1| + |v_2 - w_2|$.
By $B(v,k)$ and $\partial B(v,k)$ we denote the set of vertices $w$ for which
$|v-w|$ is at most $k$, respectively equal to $k$.
For $V, W \subset \Z^2$, we define the distance between $V$ and $W$ as
$\min\{|v-w| \, : \, v\in V, w \in W\}$.

Recall that $\calpd$ denotes the sdp distribution (that is, the distribution of the collection
$(Z_i, \, i \in \Z^2)$ defined in Subsection 1.1).
This is a distribution on
$\Om := \{0,1\}^{\Z^2}$ (with the usual $\sigma$-field).
Elements of $\Om$ are typically denoted by $\om (= (\om_i, i \in \Z^2))$, $\si$ etc.
We write $\om \leq \si$ if $\om_i \leq \si_i$ for all $i$. \\
Let $V$ be a set of vertices 
and $A$ an event. We say that $A$ lives on $V$ if 
$\om \in A$ and $\si_i = \om_i$ for all $i \in V$, implies $\si \in A$. 
And we say that $A$ is a cylinder event if $A$ lives on some finite set of vertices.
As usual, we say that $A$ is increasing if $\om \in A$ and $\om_i \leq \si_i$ for all $i$, implies 
$\si \in A$.  The first two lemma's below come from Section 2.2 and 2.4 respectively in \cite{BeB}.

\begin{lem}
Let $A$ and $B$ be two increasing cylinder events. We have
$$\calpd(A \cap B) \geq \calpd(A) \calpd(B).$$
\label{FKG}
\end{lem}
As to monotonicity, it is obvious that the sdp model has monotonicity in $\de$: If $\de_1 \geq \de_2$, then
$\calP_{p, \de_1}$ stochastically dominates
$\calP_{p, \de_2}$. Although there seems to be no `nice' monotonicity in $p$ we have the following result.

\begin{lem}
If $p_2 \geq p_1$ and $p_2 + (1-p_2) \de_2 \leq p_1 + (1-p_1)\de_1$, then
$$\calP_{p_1,\de_1} \mbox { dominates } \calP_{p_2, \de_2}.$$
\label{mon}
\end{lem}

The next result is about `almost independence' of cylinder events which live on widely separated sets.
As usual, the lattice which has the
same vertices as the square lattice but where each vertex has, besides the four
edges to its nearest neighbours, also four `diagonal edges' is called the
matching lattice (of the square
lattice). To distinguish paths and circuits in the matching lattice from those in the square lattice, we
use the terminology *-paths and *-circuits.

\begin{lem}
Let $k$ be a positive integer and let $V$ and $W$ be subsets of $\Z^2$ that have distance larger than
$2 k$. Further, let $A$ and $B$ be events which live on
$V$ and $W$ respectively. Then 

\begin{eqnarray}
& & |\calpd(A \cap B) - \calpd(A) \calpd(B)| \, \leq
\label{asindep0}  \\
& & 2 (|V| + |W|) \calpp \left(\exists \mbox{ vacant *-circuit }
\mbox{ surrounding } O \mbox{ and some vertex in } \partial B(O,k)\right).
\nonumber
\end{eqnarray}
\label{asindep}
\end{lem}

\begin{proof}
Recall how we formally defined the sdp model in terms of  random variables $X$, $Y$ and $Z$.
We use a modification of those variables: Let $X$ and $Y$ be as before, but in addition to $X^*$ and $Z$ we now
define $X^{*(f)}$ and $Z^{(f)}$ by
\begin{eqnarray}
X_i^{*(f)} &=&
\left\{ \begin{array}{ll}
& 1 \mbox{ if } X_i = 1 \mbox { and there is no $X$- occupied path from } i \mbox{ to }
\partial B(i,k) \\
& 0 \mbox{ otherwise};
\end{array}
\right. \nonumber \\
\nonumber\\
Z^{(f)}_i &=& X_i^{*(f)} \vee Y_i.
\end{eqnarray}

\noindent
Let $\calpd^{(f)}$ denote the distribution of $Z^{(f)}$.
It is clear that the random variables $Z_i^{(f)}, i \in V$ are independent of the random variables
$Z_i^{(f)}, i \in W$, and hence 
\begin{equation}
\calpd^{(f)}(A \cap B) = \calpd^{(f)}(A) \calpd^{(f)}(B).
\label{asindep1}
\end{equation}
Also note that if $Z_i \neq Z^{(f)}_i$, then the $X-$occupied cluster of $i$ intersects $\partial B(i,k)$ but 
is finite. Hence there is an $X-$vacant circuit in the matching lattice that surrounds $i$ and some site
in $\partial B(i,k)$.
Hence, since the $X-$ variables are Bernoulli random variables with parameter $p$, we have for any finite set $K$
of vertices and any event $E$ living on $K$,

\begin{eqnarray}
\label{asindep2}
& & |\calpd(E) - \calpd^{(f)}(E)| \, \leq \,
P(Z_K \neq Z^{(f)}_K) \, \leq\\  \nonumber
& & |K| \calpp \left(\exists \mbox{ a vacant  *-circuit }
 \mbox{ surrounding } O \mbox{ and some vertex in } \partial B(O,k)\right).
\end{eqnarray}

The lemma now follows easily from \eqref{asindep1} and \eqref{asindep2}
\end{proof}

%If $p > p_c$, then $1-p$, the probability that a given site is vacant, is smaller than the critical
%probability of the matching lattice. By applying the well-known exponential-decay theorems for 
%subcritical ordinary percolation to the r.h.s. of \eqref{asindep} (and noting that
%if $p<p_c$, then $\calpd(A \cap B) = \calpd(A) \calpd(B)$ for all events $A$ and $B$ that live on disjoint sets), we get

%\begin{cor}
%For all $p \neq p_c$ there is a $\psi(p) > 0$ such that for all 
%positive integers $k, l, m$ and all events 
%$A$ and $B$ that live on finite sets of vertices that have size $k$, respectively $m$, and that are separated by
%distance larger than $2 k$, we have
%$$|\calpd(A \cap B) - \calpd(A) \calpd(B)| \leq (l+m) \exp(-k \psi(p).$$
%\label{expindep}
%\end{cor}

Our last result in this section is on the uniqueness of the infinite cluster.
\begin{lem}
If $\th(p,\de) > 0$, then 
$$\calpd(\exists \mbox{ a unique infinite occupied cluster }) = 1.$$
\label{unique}
\end{lem}
\begin{proof}
From the earlier construction of the sdp model in terms of  the $X-$ and $Y$ variables, it is clear that $\calpd$ is
stationary and ergodic.
It is also clear that in the sdp model  the conditional probability that a given site is occupied given the configuration at all other sites,
is at least $\de$.
So this model has the so-called positive finite energy property. The result now follows from an extension in \cite{GaKN}
of the
well-known Burton-Keane (\cite{BuK}) uniqueness result.
%(An alternative proof, which works only for planar lattices,  is the following: 
%Since the model has positive association (Lemma \ref{fkg}) and has the symmetries of the square lattice,
% the result follows from Gandolf, Keane and ???).
\end{proof}

\label{basic}
\end{section}
\begin{section}{Partial continuity results}
In this section we first prove that in the sdp model the probabilities of cylinder events are continuous functions of
$(p,\de)$. Next we  prove that the function $\th(\cdot,\cdot)$ is continuous at $(p,\de)$ if $\th(p,\de) = 0$ or there
is an open neighborhood of $(p,\de)$ on which $\th$ is strictly positive. Note that, once we have this, the proof of
Theorem \ref{mainthm} is basically reduced to showing that if $p \neq p_c$ and $\th(p_c, \de) > 0$, then $\th(\cdot,\cdot)$
is strictly positive on an open neighborhood of $(p,\de)$. 

\begin{lem}
Let $A$ be a cylinder event. The function $(p,\de) \ra \calpd(A)$ is continuous on $[0,1]^2$.
\label{cylcont}
\end{lem}
\begin{rem}
The proof (see below) uses the well-known fact that $\th(p_c) = 0$. For many lattices (e.g. the
cubic lattice) this fact has not been proved. For those lattices the arguments below show that the function
in the statement of \ref{cylcont} is continuous on $[0,1]^2 \setminus (\{p_c\} \times [0,1])$.
\end{rem}

\begin{proof}
Let $A$ be an event which lives on some finite set $V$. Recall the construction of the sdp model in terms of
random variables $X$, $Y$ and $Z$. Let, for $\si \in \Om$, $\sigma_V$ denote the tuple $(\sigma_i, i \in V)$.
It is clear that the distribution of $X^*_V$ is a function of $p$ only, and that,
conditioned on $X^*_V$, the probability that 
$Z_V \in A$ is a polynomial (of degree $|V|$) in $\de$. Therefore it is sufficient to prove that, for each
$\al \in \{0,1\}^V$, the function $f: \,\,p \ra \calP(X^*_V = \al)$ is continuous.
Recall that the $X-$ variables are Bernoulli random variables (with parameter $p$).
Now let $0 < p_1 < p_2$. In a standard way, by introducing independent, uniformly on the interval $(0,1)$ 
distributed random variables $U_i, i \in \Z^2$, we can suitably couple two collections of Bernoulli random variables with
parameters $p_1$, respectively $p_2$. Such argument easily gives that
$|f(p_2) - f(p_1)|$ is less than or equal to the sum over $i \in V$ of 
$$\calP(U_i \in (p_1,p_2)) \, + \calP(i \mbox{ is in an infinite } p_2 \mbox{-open but not in an infinite }
p_1 \mbox{-open cluster}),$$
which equals 
$$|V| (p_2 - p_1 + \th(p_2) - \th(p_1)).$$
The lemma now follows from the continuity of $\th(.)$.
\end{proof}

\begin{prop}
Let $(p,\de) \in [0,1]^2$. If (a) or (b) below holds, the function $\th(\cdot,\cdot)$ is continuous at $(p,\de)$. 

\smallskip\noindent
(a) $\th(\cdot,\cdot) > 0$ on an open neighborhood of $(p,\de)$. \\
(b) $\th(p,\de) = 0$, 
\label{contprop}
\end{prop}

\begin{proof}
For this (and some other) results it is convenient to describe the sdp model in terms of Poisson processes:
Assign to each site, independently of the other sites, a Poisson clock with rate $1$.
These clocks govern the following time evolution:
Initially each site is vacant. Whenever the clock of a site rings, the site becomes occupied.  (If it was already occupied,
the ring is ignored). Note that if occupied sites would always remain occupied, then for each time $t$, the
configuration at time $t$ would be a
collection of independent Bernoulli random variables with parameter $1-\exp(-t)$. In particular, before and at time $t_c$,
defined by the relation $p_c = 1 - \exp(-t_c)$, there would be no infinite occupied cluster, but after $t_c$ there
would be a (unique) infinite cluster. However, we do allow occupied sites to become vacant, although only once,
as follows:
Fix a time $\tau$, a parameter of
the time evolution. At time $\tau$
all sites in the infinite occupied cluster become
vacant. (If there is no infinite occupied cluster, which is a.s. the case
if  $\tau \leq t_c$, nothing happens). After time $\tau$ we let the evolution behave as before; that is,
each vacant site becomes occupied when its Poisson clock rings.
Let, for this time evolution with parameter $\tau$,
$\hcaltat$ denote the distribution of the configuration at time $t$ ,
and let
\begin{equation}
\hth(\tau,t) = \hcaltat(O \mbox{ is in an infinite occupied cluster }).
\label{htheq}
\end{equation}
It is easy to see that
\begin{equation}
\hcaltat = \calpd,
\label{tatpdeq}
\end{equation}
where 
$p = 1 - \exp(-\tau)$ and $\de = 1 -\exp(-(t-\tau))$.
It is also easy to see that  $\hcaltat$ is stochastically decreasing in $\tau$ and stochastically increasing in $t$.
In fact this is the key behind Lemma \ref{mon}.

Now we come back to the proof of Proposition \ref{contprop}.
From \eqref{htheq} and \eqref{tatpdeq} we get (since the map between pairs $(p,\de)$ and $(\tau,t)$ in \eqref{tatpdeq}
is continuous) that this proposition is equivalent to saying that if $\hth(\tau,t) = 0$ or
$\tau \neq t_c$ and $\hth$ is strictly positive on
an open neighborhood of $(\tau,t_c)$, then $\hth$ is continuous at $(\tau,t)$.
To prove this equivalent form of Proposition \ref{contprop} we use ideas from \cite{BeK}.
The introduction of pairs $(\tau,t)$ as replacement of $(p,\de)$ not only has the advantage that, as we already
saw, we now have a more suitable form of monotonicity, but, more importantly, that we now have a more 'detailed'
structure (the Poisson processes) in the background which gives the appropriate `room' needed to get a 
suitable modification of the arguments in \cite{BeK}.

Let $(\tau,t)$ be as above. Divide the parameter space in four `quadrants', numbered $I$ to $IV$:

\begin{eqnarray}
\nonumber
& & I  := [0,\tau] \times [t,\infty), \\ \nonumber
& & II := [\tau,\infty) \times [t,\infty), \\ \nonumber
& & III:= [\tau,\infty) \times [0,t], \\ \nonumber
& & IV := [0,\tau] \times [0,t]. \nonumber
\end{eqnarray}
Note that it is sufficient to prove that for each monotone sequence $(\tau_i, t_i)_{i \geq 0}$ that lies
in one of the above quadrants and converges to $(\tau,t)$, one has 
$$\lim_{i \ra \infty} \hth(\tau_i,t_i) = \hth(\tau,t).$$
We handle each of the quadrants separately.

\smallskip\noindent
Quadrant I) This is easy and corresponds to the (easy) proof of right continuity
of ordinary percolation: Let $(\tau_i)$ be a monotone sequence which converges from below to $\tau$ and let $(t_i)$ be 
a monotone sequence
which converges from above to $t$. 
Let $A_n$ denote the event that there is an occupied path from $O$ to
$\partial B(O,n)$.
By monotonicity and Lemma \ref{cylcont} we have that
\begin{eqnarray}
& & \mbox{ For each } i, \hcaltaiti(A_n) \downarrow \hth(\tau_i,t_i) \mbox{ as } n \rightarrow \infty; \\
& & \hcaltat(A_n) \downarrow \hth(\tau,t) \mbox{ as } n \rightarrow \infty; \\
& & \mbox{ For each } n, \hcaltaiti(A_n) \downarrow \hcaltat(A_n) \mbox{ as } i \rightarrow \infty,
\label{I2}
\end{eqnarray}
From these three statements it is easy to see that $\hth(\tau_i, t_i)$ tends to $\hth(\tau,t)$ as $i \ra \infty$.

Quadrant III)
Let $(\tau_i)$ be a monotone sequence which converges from above to $\tau$ and $(t_i)$ a
monotone sequence
which converges from below to $t$. 
By the earlier monotonicity arguments, the sequence  $\hth(\tau_i,t_i)$ is increasing in $i$, and has
a limit smaller than or equal to $\hth(\tau,t)$.
So for the situation where $\hth(\tau,t) = 0$, the proof is done.
Now we handle the other situation: we assume $\hth$ is positive in an open neighborhood of $(\tau,t)$. For this
situation considerable work has to be done. 
Note that in the dynamic description given earlier in this section, the underlying Poisson
processes were the same for each choice of the model parameter $\tau$. This allows us (and we
already used this to derive some monotonicity properties) to couple the models with the different
$\tau_i$'s and $\tau$. 

Let, for $s < u$, $C_{s,u}$ denote the occupied cluster of site $O$ at time $u$ in the process with parameter $s$
(that is, under the time evolution where the infinite occupied cluster is destroyed at time $s$).
Further, we use the notation $\om(s,u)$ for the configuration at time $u$ in that model.
It is also convenient to consider $\om(u)$, the configuration at time $u$ in the model where {\it no}
destruction takes place. (So, $\om_v(u), v \in \Z^2$, are independent $0-1$ valued random variables, each being
$1$ with probability $1- \exp(-u)$).
Again we emphasize that all these models are defined in terms of the same
Poisson processes.
From monotonicity (note that $C_{\tau_i, t_i} \subset C_{\tau_{i+1}, t_{i+1}}$ for all $i$) it is clear that 
$$\lim_{i \ra \infty} \hth(\tau_i, t_i) = \calP(\exists \, i \, |C_{\tau_i,t_i}| = \infty),$$
and
\begin{equation}
\hth(\tau,t) - \lim_{i \ra \infty} \hth(\tau_i, t_i) = \calP(|C_{\tau,t}| = \infty, \, \forall i \,
|C_{\tau_i,t_i}| < \infty).
\label{diff}
\end{equation}
So we have to show that the r.h.s. of \eqref{diff} is $0$.
Fix a $j$ with the property that $\hth(\tau_j,t_j) >0$.
Such $j$ exists by the condition we assumed for $(\tau,t)$.
To show that the r.h.s. of \eqref{diff} is $0$, it is sufficient (and necessary) to prove the following
claim:

\smallskip\noindent
{\bf Claim}\\
{\it Apart from an event
of probability $0$, the event $\{|C_{\tau,t}| = \infty\}$ is contained in the event that
there is a $k > j$ for which $|C_{\tau_k,t_k}| = \infty$.}

\smallskip\noindent
So suppose $|C_{\tau,t}| = \infty$.
By our choice of $j$ we may assume that $\om(\tau_j,t_j)$ has an infinite occupied cluster, and by Lemma \ref{unique}
that this cluster is unique. We denote it by $I_j$. If $O \in I_j$ we are done.
From monotonicity and the uniqueness of the infinite cluster (see Lemma \ref{unique}), we have $I_j \subset C_{\tau,t}$.
Hence there is a finite path $\pi$ from $O$ to some site in $I_j$ such that $\om(\tau,t) \equiv 1$ on $\pi$.
Since, a.s. there are no vertices whose clock rings exactly at time $t$ or $\tau$,
we may assume that for every site $v$ on $\pi$, (a) or (b) below holds:

\smallskip\noindent
(a) The clock of $v$ rings in the interval $(\tau,t)$. \\
(b) $\om_v(\tau) = 1$ but the occupied cluster of $v$ in $\om(\tau)$ is finite.

\smallskip\noindent
If (a) occurs we define:

$$i_v := \min\{i \, :\, i \geq j \mbox{ and the clock of } v \mbox{ rings in } (\tau_i,t_i)\}.$$

\noindent
Note that then, by the monotonicity of the sequence $(\tau_i,t_i)$, the clock of $v$ rings in
the interval $(\tau_l, t_l)$ for all $l \geq i_v$.
If (a) does not occur, (b) occurs, and hence there is a finite set $K_v$ of sites on which $\om(\tau) = 0$ and which
separates $v$ from $\infty$. Then we use the following alternative definition of $i_v$:

$$i_v := \min\{i \, : \, i \geq j \mbox{ and } \om(\tau_i) \equiv 0 \mbox{ on } K_v\}.$$
This minimum exists since $K_v$ is finite and (again) we assume that no Poisson clock rings exactly at 
time $\tau$.
Now let 

$$k := \max_{v \in \pi} i_v,$$
which exists since $\pi$ is finite.

From the above procedure it is clear that $\om_v(\tau_k, t_k) = 1$ for all $v$ on $\pi$. Further, since $k \geq j$ and
by monotonicity,  also $\om_v(\tau_k, t_k) = 1$ for all $v \in I_j$. Since $\pi$ is a path from
$O$ to $I_j$ this implies that $I_j$ is contained in $C_{\tau_k, t_k}$ and hence that $|C_{\tau_k, t_k}| = \infty$.
This proves the Claim above.

\smallskip\noindent
Quadrants II) and IV) \\
The required results for these quadrants follow very easily from monotonicity and the above results for quadrants
I and III: Let $(\tau_i,t_i)$ be a sequence in quadrant II that converges to $(\tau,t)$.
We have, by earlier stated monotonicity properties,

$$\hth(\tau_i,t) \leq \hth(\tau_i,t_i) \leq \hth(\tau,t_i).$$

Since the sequence $(\tau_i,t))$ lies in quadrant III and the sequence $(\tau,t_i)$ lies in quadrant I,
the upper and lower bound both converge to $\hth(\tau,t)$. This completes the treatment of quadrant II.
Quadrant IV is treated in the same way.
This completes the proof of Proposition \ref{contprop}
\end{proof}
\end{section}
\begin{section}{An RSW-type result}
For our main result we need to prove that if the crossing probability of an $n$
by $n$ square  goes to $1$ as $n \ra \infty$,
then also the crossing probability of a (say) $3 n \times n$ rectangle
in the `difficult direction' goes to $1$ as $n \ra \infty$.
Such (and stronger) results were proved for ordinary percolation in the late nineteen seventies by Russo, and by
Seymour and Welsh, and therefore became known as RSW theorems.
%In fact, Russo, and Seymour and Welsh
%proved much stronger statements. 
%For instance, they showed that there is a strictly increasing function $f : [0,1] \times [0,1]$ such that for
%each $n$ the probability to cross a $3n \times n$ reactangle is larger than or equal to 
%$$f(\mbox{ the probability of crossing an } n\times n \mbox{ rectangle }).$$
Their proofs used careful conditioning on
the lowest horizontal crossing in a rectangle, after which the area
above that crossing was treated, and a new, vertical crossing in that area was `constructed'.
Such arguments work
for ordinary percolation because there the above mentioned area can be treated as `fresh' territory. However, 
they usually break down in situations where we have dependencies, as in the sdp model. 

Recently, Bollob\'{a}s and Riordan (\cite{BoR}) made significant progress on these matters. For the so-called Voronoi 
percolation model they proved an RSW type result. That result is one of the main ingredients in their proof
that the critical probability for Voronoi percolation equals $1/2$ (which had been conjectured but stayed open
for a long time). Although they explicitly proved their RSW type result only for the Voronoi model, their proof
works (as they remark in their paper) for a large class of models.
The result we needed is a little stronger than that of \cite{BoR}. The rest of this section is organised as follows.
First we give a short introduction to Voronoi percolation. Then we state the above mentioned 
RSW-like theorem of \cite{BoR}, and point out where and how its proof needs to be modified to obtain
the stronger version. Finally we state the analog for the sdp model and explain why the proof for the Voronoi model
works for this model as well. 
\begin{subsection}{The Voronoi percolation model}
We start with a brief description of the Voronoi percolation model. 
The (random) Voronoi percolation model is as follows: Let $Z$ denote the (random) set of points in a Poisson point
process with density $1$ in the plane. This set gives rise to a random Voronoi tessellation of the plane:
Assign to each $z \in Z$ the set of all $x \in \R^2$ for which $z$ is the nearest 
\medskip\noindent point in $Z$. The closure of this set is called the (Voronoi) cell of $z$. 
It is known that (with probability $1$) each Voronoi cell is a convex polygon, and that two cells are either
disjoint or share an entire edge. In the latter case the two cells are said to be neighbours or adjacent.
This notion of adjacency gives, in a natural way, rise to the notion of paths, clusters etc.

Now consider the percolation model where each cell, independently of everything else, is
coloured black with probability $p$ and white with probability $1-p$. Based on analogies with ordinary percolation
(in particular with the self-matching property of the usual triangular lattice) it has been conjectured for
a long time that the critical value for this percolation model is $1/2$: for $p < 1/2$ there is (a.s.) no
infinite black cluster, but for $p > 1/2$ there is an infinite black cluster (a.s.). As we said before,
this was recently proved rigorously
by Bollob\'{a}s and Riordan (\cite{BoR}), and a key ingredient in
their proof
is an ingenious RSW-like result.
%In the following subsection we discuss this result and explain how, with only
%a little extra work, a somewhat stronger form can be obtained. Finally we point out that a similar result holds
%for the sdp model. 
\end{subsection}
\begin{subsection}{The RSW-like result for Voronoi percolation}
As in \cite{BoR} we define, for the Voronoi percolation model with parameter $p$,
$f_p(\rho,s)$ as the probability that there is a horizontal black crossing
of the rectangle
$[0,\rho s] \times [0,s]$.  The following is Theorem 12 in \cite{BoR}

\begin{thm} (Bollob\'{a}s and Riordan)
Let $0<p<1$ be fixed. 

\begin{eqnarray}
\mbox{ If } & & \liminf_{s \ra \infty} f_p(1,s) > 0, 
\label{brcondition}
\\ \nonumber
\mbox{ then } & & \limsup_{s \ra \infty} f_p(\rho,s) > 0 \mbox{ for all } \rho >0.
\end{eqnarray}
\label{brthm}
\end{thm}

Studying the proof we realised that the condition can be weakened, so that the following theorem is obtained:

\begin{thm}
Let $0 < p <1$ be fixed.  
\begin{eqnarray}
\nonumber
\mbox{ If } & & \limsup_{s \ra \infty} f_p(\rho,s) > 0 \mbox { for {\it some} } \rho > 0, \\ \nonumber
\mbox{ then }& & \limsup_{s \ra \infty} f_p(\rho,s) > 0 \mbox{ for all } \rho >0. \nonumber
\end{eqnarray}
\label{strongbrthm}
\end{thm}

%In fact, the conclusion, as formulated in the BollR paper is for $\rho > 1$ instead of $\rho >0$, but that is not
%a restriction since $f(\rho,s)$ is strictly non-increasing in $\rho$. 
%We formulated their result in the
%form above (Theorem \ref{brthm}), to make it easier to compare with the stronger Theorem \ref{strongbrthm}. 
%The differences between the two theorems are that the second has a $\limsup$ instead of $\liminf$
%in its condition, and that it allows $ \rho < 1$ in the condition.

\smallskip\noindent
As we shall point out, this somewhat stronger Theorem \ref{strongbrthm} can be proved in almost the same way as
Theorem \ref{brthm}. But see Remark \ref{remdirect} about the global structure of the proof.
% is by contradiction. In our personal view, one can see more
%clearly that minor changes give the stronger version when instead of a proof by contradiction we give a proof with a more
%`direct' structure, which will be done below.
%\label{remdiff}
%\end{rem}
First note that Theorem \ref{strongbrthm} is (trivially) equivalent to 
the following:

\begin{thm}
Let $0 < p <1$ be fixed.  

\begin{equation}
\mbox{ If }\limsup_{s \ra \infty} f_p(\rho,s) = 0 \mbox { for {\it some} } \rho > 0,
\label{condsbrthm}
\end{equation}
then
\begin{equation}
\limsup_{s \ra \infty} f_p(\rho,s) = 0 \mbox{ for all } \rho >0.
\label{conclusion}
\end{equation}
\label{eqsbrthm}
\end{thm}

This is the form we will prove, following (with some small changes) the steps in \cite{BoR}.

\begin{proof} (Theorem \ref{strongbrthm} and \ref{eqsbrthm}).
Since $p$ is fixed we will omit it from our notation. In particular we will write $f$ instead of $f_p$.

First we will rewrite the condition \eqref{condsbrthm} in Theorem \ref{eqsbrthm}:
If $\limsup_{s\ra\infty} f(\rho,s) = 0$ for some $\rho \leq 1$ then, since $f(\rho,s)$ is decreasing in $\rho$,
this $\limsup$ is $0$ for all $\rho > 1$. 
Moreover,
the well-known pasting techniques from ordinary percolation show easily that if 
$\limsup_{s \ra \infty} f(\rho,s) > 0$ for some $\rho > 1$, then this $\limsup$ is positive for all $\rho' > \rho$,
and hence (using again monotonicity of $f$ in $\rho$) for all $\rho >1$.
Equivalently, if $\limsup_{s \ra \infty} f(\rho,s) = 0$ for some $\rho > 1$, then this limit equals $0$ for all $\rho >1$.
Hence, the condition in Theorem \ref{eqsbrthm} is equivalent to

\begin{equation}
\limsup_{s \ra \infty} f_p(\rho,s) = 0 \mbox { for all } \rho > 1,
\label{eqcondition}
\end{equation}
and this is also equivalent to condition (3) in Section 4 of \cite{BoR}:

\begin{equation}
\limsup_{s \ra \infty} f_p(\rho,s) = 0 \mbox { for some } \rho > 1,
\label{bollcond3}
\end{equation}

We will assume \eqref{eqcondition} (or its equivalent form \eqref{bollcond3})
and show how, following basically the proof of Theorem \ref{brthm}, the equation in \eqref{conclusion} can be derived from it
for all $\rho > 1/2$.
%In fact, we will point out that following  the proof in \cite{BoR}
% step by step, with a few small changes, shows that  \eqref{conclusion} 
%holds for all $\rho > 1/2$.
Then we make clear that, for each $k$, very similar arguments work for $1/k$
instead of $1/2$, which completes the proof of Theorem \ref{eqsbrthm}.

\begin{rem}
Bollob\'{a}s and Riordan prove their theorem by contradiction:
They assume (as we do here) \eqref{bollcond3} above,
and, moreover they assume 
\eqref{brcondition} (equation (2) in Section 4 of their paper).
Then, after a number of steps (claims), they reach a contradiction, which completes the proof.
However, most of these steps do not use the `additional' assumption \eqref{brcondition} at all.
We found a `direct' (that is, not by contradiction) proof, as sketched below, more clarifying
since it leads more easily to further improvements. For our goal most of the steps (claims) in the proof in \cite{BoR}
remain practically unchanged. Therefore we (re)write only some of them in more detail (Claim 1 is stated
to give an impression of the start of the proof, and Claim 4 because that already gives a good indication of the
strong consequences of \eqref{bollcond3}). For the other claims we only describe which
changes have to be made for our purpose.
\label{remdirect}
\end{rem}

First some notation and terminology: $T_s$ is defined as the strip $[0,s] \times \R$.
An event is said to hold {\it with high probability}, abbreviated {\it whp} if its probability goes to $1$ as
$s \ra \infty$ (and all other parameters, e.g. $p$ and $\epsilon$ are fixed).

\smallskip\noindent
{\bf Claim 1} (Claim 12.1 in \cite{BoR}). \\
{\it
Let $\ep > 0$ be fixed, and let $L$ be the line-segment $\{0\} \times [-\ep s, \ep s]$.   
Assuming that \eqref{bollcond3} holds, the probability that there is a black path $P$ in $T_s$
starting from $L$ and going outside $S' = [0,s] \times [-(1/2 + 2 \ep) s, (1/2 + 2 \ep)s]$
tends to $0$ as $s \ra \infty$.
} \\

\smallskip\noindent
This claim is exactly the same as in \cite{BoR}, except that in their formulation not only \eqref{bollcond3} but
also \eqref{brcondition} is assumed. However, their proof of this claim does not use the latter assumption.

The above, quite innocent looking claim, leads step by step to stronger and eventually very strong
claims. We will not rewrite {\bf Claim 2} and {\bf Claim 3}; like Claim 1, they are exactly the same as their corresponding
Claims (12.2 and 12.3 respectively) in \cite{BoR}, except that we do not assume
\eqref{brcondition}. And, again, the proof remains as in \cite{BoR}.

\smallskip\noindent
{\bf Claim 4} (Claim 12.4 in \cite{BoR})
{\it
Let $C > 0$ be fixed, and let $R=R_s$ be the $s$ by $2 C s$ rectangle $[0,s] \times [-C s, C s]$.
For $0 \leq j \leq 4$, set $R_j = [j s/100, (j + 96)s/100] \times [-C s, C s]$.
Assuming that \eqref{bollcond3} holds, {\it whp} every black path $P$ crossing $R$ horizontally
contains $16$ disjoint black paths $P_i, \, 1 \leq i \leq 16$, where each $P_i$ crosses some $R_j$
horizontally.
}

\smallskip\noindent
Again, in the formulation in \cite{BoR} also \eqref{brcondition} is assumed, but this is not used in the proof.
Following \cite{BoR} we now define, for a rectangle $R$, the random variable $L(R)$ as the minimum length of a black path
crossing $R$ horizontally. (More precisely, it is the minimum length of a piecewise-linear black curve that crosses
$R$ horizontally). If there is no horizontal black crossing of $R$ we take $L(R) = \infty$. 
A complicating property of $L$ is that if $R_1$ and $R_2$ are two disjoint rectangles, $L(R_1)$ and $L(R_2)$
are not independent (no matter how large the distance between the two
rectangles).
Therefore, below Claim 12.4 in their paper, Bollob\'{a}s and Riordan introduce a suitable modification $\tL$.
The key idea is that whp the colours inside a rectangle with length and width of order $s$,
are completely determined by the Poisson points within distance of order $o(s)$ of the rectangle.

There are many suitable choices of $\tL$, and we will not rewrite the precise definition given in
\cite{BoR}, but only highlight the following three key properties (which neither use \eqref{bollcond3} nor \eqref{brcondition}):
%for any $s$ by $2 s$ reactangle $R_s$, the event $E_{\mbox{dense}}(R_s)$ 
%that for every $x \in \R^2$ at distance smaller than $2 \sqrt{\log s}$ from $R_s$, there is a Poisson point 
%at distance smaller than $2 \sqrt{\log s}$ from $x$. It is easy to see that $E_{\mbox{dense}}(R_s)$ holds whp.
%
%The modification $\tL(R_s)$ is now defined as follows: If the event $E_{\mbox{dense}}(R_s)$ holds,
%$\tL(R_s) = L(R_s)$, otherwise it is set to $0$.
%It is easy to see that whp  $E_{\mbox{dense}}(R_s)$ holds, and hence that 

\begin{equation}
\tL(R_s) = L(R_s), \,\, \mbox{whp},
\label{Leq}
\end{equation}

\begin{equation}
\tL(R_s) \geq s, \,\, \mbox{whp},
\label{wL>s}
\end{equation}
and:

\smallskip\noindent
{\bf Claim 5.} (Claim 12.5 in \cite{BoR}). {\it Let $R_1$ and $R_2$ be two $s$ by $2 s$ rectangles, separated by a distance
of at least $s/100$. If $s$ is large enough, then the random variables $\tL(R_1)$ and $\tL(R_2)$ are 
independent.}

\begin{rem}
In fact, the independence property of $\tL$ is only used in the proof of Claim 12.6 in \cite{BoR}, and
there it could be replaced by the following property (which follows from \eqref{Leq} and Claim 5):

\smallskip\noindent
For each $\varepsilon > 0$ there is a $u$ such that for all $s > u$ and all $s$ by $2 s$ rectangles
$R_1$ and $R_2$ that are separated by a distance
of at least $s/100$, we have 
$$ \sup_{x ,y >0} | P(L(R_1) < x, \, L(R_2) < y) \, - \, P(L(R_1) < x) \, P(L(R_2) < y) | < \varepsilon.$$
\label{asindeprem}
\end{rem}
 
\smallskip\noindent
Now choose an arbitrary number $\heta$ smaller than $10^{-4}$.
This deviates from the choice of $\eta$ by Bollob\'{a}s and Riordan, who add an extra condition, related to their assumption 
of \eqref{brcondition}. Define
%Recall that we want to avoid the use of that assumption.
%Next define

\begin{equation}
\hat{t}(s) = \sup\{s \, : \, \calP(\tL(R_s) < x) \leq \heta\}.
\label{tseq}
\end{equation}

This definition of $\hat{t}$ is the same in form as that of $t$ in \cite{BoR}
(see two lines below (16) in \cite{BoR}); however our way of choosing $\heta$ was different.
A consequence of this difference is that, in our setup, $\hat{t}(s)$ can be $\infty$.
%(namely, if the probability
%of the event that $tL(R_s) = L(R_s)$
%$E_{\mbox{dense}}(R_s)$ holds
%and there is no horizontal black crossing of $R_s$ is larger than $1 - \heta$.
As in \cite{BoR}, we do have that

\begin{equation}
\hat{t}(s) \geq s, \mbox{ for all }s.
\label{bound-s}
\end{equation}

{\bf Claim 6.} (Claim 12.6 in \cite{BoR}). \\
{\it Let $R_s$ be a fixed 0.96 s by $2 s$ rectangle. If \eqref{bollcond3} holds, then}

$$\calP\left(L(R_s) < \hat{t}(0.47 s)\right) \leq 200 \heta^2.$$

This statement is the same as in \cite{BoR}, except that in \cite{BoR} also \eqref{brcondition} is assumed, and that we use
$\hat{t}$ and $\heta$ instead of $t$, respectively $\eta$. The proof is the same as in \cite{BoR}.

From the above (in particular Claim 4, Claim 6 and \eqref{Leq}), the following quite startling Proposition
(which, essentially is equation (18) in \cite{BoR}) now follows quite easily.

\begin{prop}
({\em Corresponds with (18) in \cite{BoR}}).
If \eqref{bollcond3} holds, then, for all sufficiently large $s$,

\begin{equation}
\hat{t}(s) \geq 16 \hat{t}(0.47 s).
\label{trecursion}
\end{equation}

\label{t-rec-prop}
\end{prop}

\begin{proof}
Practically the same as the proof of equation (18) in \cite{BoR}.
\end{proof}

Now Theorem \ref{eqsbrthm} follows in a few lines from this proposition and the definition of $\hat{t}(s)$:
It is easy to show (see the arguments below (18) in \cite{BoR}) that \eqref{bound-s} and Proposition \ref{t-rec-prop}
together imply that
$\hat{t}(s) > s^3$, for all sufficiently large $s$.
Hence, by the definition of $\hat{t}(s)$ (and by \eqref{Leq}) we get that,
for all sufficiently large $s$, 
\begin{equation}
\calP(L(R_s) < s^3) \leq \hat{\eta}.
\label{s3bound}
\end{equation}

It is also easy to show  (see the arguments below equation (19) in \cite{BoR}) that

\begin{equation}
\calP(s^3 \leq L(R_s) < \infty) \ra 0 \mbox{ as } s \ra \infty.
\label{s3bound2}
\end{equation}

Combining \eqref{s3bound} and \eqref{s3bound2}, and recalling that $L(R_s) = \infty$ iff there is no
horizontal black crossing of $R_s$, immediately gives

$$\limsup_{s \ra \infty} \calP(\exists \mbox{ a horizontal black crossing of } R_s) \leq \hat{\eta}.$$

Now, since $\eta$ was an arbitrary number between $0$ and $10^{-4}$, we get \\
$\lim_{s \ra \infty} \calP(\exists \mbox{ horizontal black crossing of } R_s) = 0$,
that is,

\begin{equation}
f(1/2, s) \ra 0, \mbox{ as } s \ra \infty.
\label{thm-for-2}
\end{equation}

Note that in the last part of the above arguments (after Claim 4) we worked in particular with $s$ by $2 s$ rectangles.
A careful look at the arguments shows that the choice of this factor $2$ is, in fact, immaterial: if we would take
$s$ by $3s$ rectangles or, more generally, fix an $N \geq 2$ and take $s$ by $N s$ rectangles, the arguments
remain practically the same. To see this, one can easily check that in Claims 1 - 4 (Claims 12.1 - 12.4 in \cite{BoR})
the factor $2$ plays no role at all:
here the rectangles under consideration are $s$ by $2 C s$, where $C$ is a fixed but arbitrary positive number.
Further, the proof of Claim 5 remains the same when, for some fixed positive number $C$, we replace the factor $2$
by $2 C$.
And, the  definition of $\hat{t}(s)$ (see \eqref{tseq}), which was given in terms of $s$ by $2 s$ rectangles, has,
for each $C >0$, an obvious analog for $s$ by $2 C s$ rectangles:

\begin{equation}
t_C(s) :=  \sup\{s \, : \, \calP(\tL(R^C_s) < x) \leq \heta\},
\label{tCseq}
\end{equation}
where, for each $s$, $R^C_s$ is some fixed $s$ by $2 C s$ reactangle.

In the generalization of Claim 6 (Claim 12.6 in \cite{BoR}) we now fix $C \geq 1$, and take
$R_s^C:= [0, 0.96 s] \times [-C s, C s].$ 
In the proof
of this Claim we have to replace, on the vertical scale, $s$ by $C s$.
For instance, the segments $L_i$, which in the original proof in \cite{BoR} have length $0.02 s$, will now have length
$0.02 C s$, and $R_0$ and $R_1$ which in the original proof are $0.47 s$ by $2 \times 0.47 s$
rectangles, are now $0.47 s$ by $2 C 0.47 s$ rectangles.
In this way we get if \eqref{bollcond3} holds,  for each fixed $C > 1$ the following analog of \eqref{thm-for-2}:

\begin{equation}
f(1/(2 C), s) \ra 0, \mbox{ as } s \ra \infty.
\label{thm-for-2C}
\end{equation}

This proves 
Theorem \ref{eqsbrthm} and hence Theorem \ref{strongbrthm}.

\end{proof}

In the above we were dealing with black horizontal crossings. Obviously, completely analogous results hold for
white horizontal crossings:
If we denote (for a fixed value of the parameter $p$ of the Voronoi percolation model), the probability of a
vertical white crossing of a given $\rho s$ by $s$ rectangle by
$g(\rho,s)$, we have that if $\lim_{s \ra \infty} g(\rho,s) = 0$ for {\it some} $\rho >0$, then
this limit is $0$ for {\it all} $\rho>0$. Since a rectangle has either a horizontal
black crossing or a vertical white crossing (and hence
$g(\rho,s) = 1 - f(\rho, s)$) this gives:

\begin{cor}
\begin{eqnarray}
\mbox{ If } & & \lim_{s\ra \infty} f(\rho,s) = 1 \mbox{ for some } \rho > 0,\\ \nonumber
\mbox{ then } & & \lim_{s\ra \infty} f(\rho,s) = 1 \mbox{ for all } \rho > 0, \nonumber
\end{eqnarray}
\label{cor-sbrthm}
\end{cor}

\end{subsection}
\begin{subsection}{An RSW analog for self-destructive percolation}
In the previous subsection we considered (and somewhat strengthened) an RSW-like result of Bollob\'{a}s and
Riordan (\cite{BoR}) for the Voronoi percolation model. Only a
few properties of the model are used in its proof. As remarked in \cite{BoR} (at the end of Section 4; see also 
\cite{BoR2}, Section 5.1),
these properties are basically the
following: First of all, crossings of rectangles are defined in terms of `geometric paths' in such a way that,
for example, horizontal and vertical black crossings meet, which enables to form longer paths by pasting together
several small paths. Further, a form of FKG is used (e.g. that events of the form `there is a black path from
$A$ to $B$' are positively correlated. Also some symmetry is needed. Bollob\'{a}s and Riordan say that "invariance of the
model under the symmetries of $\Z^2$ suffices, as we need only consider rectangles with integer coordinates". 
Finally, some form of asymptotic independence is needed
(see Remark \ref{asindeprem}).
%It is certainly
%enough that there is a $\la >0$ such that for any events $A$ and $B$ which depend on the colours in a finite region
%$V_A$ respectively $V_B$, we have 
%$$|P(A \cap B) - P(A) P(B)| \leq \exp(-\la d(V_A,V_B),$$
%where $d(V_A,V_B)$ is the distance between $V_A$ and $V_B$.
Similar considerations hold wrt the somewhat stronger Theorem \ref{strongbrthm}. \\
Using the results in Section \ref{basic}, is not difficult to see that the sdp model has the above mentioned properties: 

\begin{itemize}
\item The indicated geometric properties are just the well-known intersection properties of paths in
the square lattice (and in its matching lattice). 

\item Lemma \ref{FKG} gives the needed FKG-like properties.

\item Asymptotic independence: 
Note that for $p \leq p_c$ the sdp model is an ordinary percolation model, where
this property is trivial. If $p>p_c$, then $1 - p$ is smaller than the critical probability of
the matching lattice. In that case the needed asymptotic independence (of the form described in Lemma 4.5) comes
from Lemma \ref{asindep} and the well-known
exponential decay theorems for ordinary subcritical percolation. 

\item The sdp model on the square lattice clearly has all the symmetries of $\Z^2$.
\end{itemize}
%COMMENT: SHOULD WE GIVE MORE DETAILS?

Further, to carry out for the sdp model the analog of the arguments that led from
Theorem \ref{eqsbrthm} to Corollary \ref{cor-sbrthm}, we note 
that the random collection of vacant sites on the matching lattice clearly also has the above mentioned properties.
So we get the following theorem for the sdp model:

\begin{thm}
The analogs of Theorems \ref{strongbrthm} and \ref{eqsbrthm} and Corollary \ref{cor-sbrthm} hold for the self-destructive
percolation model.
In particular, let
for the sdp model with parameters $p$ and $\de$, 
$f(\rho,s)= f_{p,\de}(\rho,s)$ denote
the
probability that there is an occupied horizontal crossing
of a given $\rho s \times s$ rectangle. We have

%\begin{eqnarray}
%\mbox{ If } & & \limsup_{s \ra \infty} f(\rho,s) > 0 \mbox { for {\t some} } \rho > 0, \\ \nonumber
%\mbox{ then } & &  \limsup_{s \ra \infty} f(\rho,s) > 0 \mbox{ for all } \rho >0. \nonumber
%\end{eqnarray}
%\label{strongbrthm-sdp}
%\end{thm}
\begin{eqnarray}
\mbox{ If } & & \lim_{s\ra \infty} f(\rho,s) = 1 \mbox{ for some } \rho > 0, \\ \nonumber
then & & \lim_{s\ra \infty} f(\rho,s) = 1 \mbox{ for all } \rho > 0. \nonumber
\end{eqnarray}
\label{thm-rsw-sdp}
\end{thm}

In the next sections this result will play an important role in the completion of the proof of
Theorem \ref{mainthm}. In particular, in Section \ref{fscriterion} it will be used to prove a finite-size criterion
for supercriticality of the sdp model.
\end{subsection}
\end{section}

\begin{section}{A finite-size criterion}
The main result of this section is a suitable finite-size criterion for supercriticality of the sdp model.
The overall structure of the argument is similar to that in ordinary percolation (see \cite{CC}), but
the dependencies in the model require extra attention. One of the main ingredients, a suitable RSW-like theorem
for this model, was obtained in the previous section.

\begin{thm}
Let $f = f_{p,\de}$ as in Theorem \ref{thm-rsw-sdp}.
There is a universal constant $\al >0$ and there is a decreasing function 
$\hN \, : \, (p_c,1) \ra \N$ such that for all $p > p_c$ and all $\de >0$
the following two assertions, (i) and (ii) below, are equivalent.

\begin{eqnarray}
& & \mbox{i. } \,\, \th(p,\de) > 0. \\ \nonumber
& & \mbox{ii.} \exists n \geq \hN(p) \mbox{ such that } f_{p,\de}(3,n) > 1-\al. \nonumber
\end{eqnarray}

\label{fsizethm}
\end{thm}
\begin{rem}
In ordinary percolation $\hN(p)$ can be taken constant $1$.
Remark \ref{remfinal} below explains the impact of this difference.
\end{rem}

\begin{proof}
Consider for each $n \in \N$ the events

\begin{eqnarray}
& & A = \{\exists \mbox{ a vertical vacant *-crossing of } [0,9 n] \times [0,3 n]\}; \\ \nonumber
& & B = \{\exists \mbox{ a vertical vacant *-crossing of } [0,9 n] \times [0, n]\}; \\ \nonumber
& & C = \{\exists \mbox{ a vertical vacant *-crossing of } [0,9 n] \times [2 n, 3 n]\}. \\ \nonumber
\end{eqnarray}

Let $h(\rho,n)$ denote the probability of a vertical vacant crossing (in the matching lattice) of 
a $\rho \, n$ by $n$ box. So, $h(\rho,n) = 1 - f(\rho,n)$.
Clearly, 
$\calpd(B) = \calpd(C) = h_{p,\de}(9,n)$ and $\calpd(A) = h_{p,\de}(3,3 n)$. It is also clear that $A \subset B \cap C$.
From this, Lemma \ref{asindep},
the fact that the r.h.s. of \eqref{asindep0} is decreasing, and
the well-known exponential decay results for ordinary subcritical
percolation applied to \eqref{asindep0}, it follows that there is an increasing,
function $\phi  \, : \,  (p_c,1) \ra (0, \infty)$ such that for all $p > p_c,$  

\begin{equation}
h(3, 3 n) \leq h(9,n)^2 + \exp(- n \phi(p)).
\label{quadineq}
\end{equation}

Further note that if the event $B$ occurs, there must be a vacant vertical *- crossing of one of the rectangles
$[0,3 n] \times [0,n]$,
$[ 2n, 5 n] \times [0,n]$, $[4 n, 7 n] \times [0, n]$, $[6 n, 9 n] \times [0,n]$,
or a vacant horizontal *-crossing of one of the rectangles
$[2 n,3 n] \times [0,n]$,
$[4 n,5 n] \times [0,n]$,
$[6 n,7 n] \times [0,n]$. 

Hence
\begin{equation}
h(9,n) \leq 4 h(3,n) + 3 h(1,n) \leq 7 h(3,n),
\label{sumineq}
\end{equation}
which combined with \eqref{quadineq} gives
\begin{equation}
h(3, 3 n) \leq 49 h(3,n)^2 + \exp(-n \phi(p)). 
\label{recineq1}
\end{equation}

Take $\al$ so small that $49 \al^2 < \al/4$. Let, for each $p > p_c$, 
$\hN(p)$ be the smallest positive integer for which
$$\exp(-\hN(p) \phi(p)) < \al/4.$$
Sine $\phi$ is increasing, $\hN$ is decreasing in $p$.

\smallskip\noindent
Now suppose $p > p_c$ and $\de \in (0,1)$ are given and suppose that (ii) holds. So there exists an $n$ that
satisfies:
\begin{equation}
\exp(-n \phi(p)) < \al/4 \mbox{ and } h(3,n) < \al.
\label{conditionii}
\end{equation}

From \eqref{recineq1}, \eqref{conditionii} and the choice of $\al$ we get

\begin{equation}
h(3, 3 n) \leq 49 \al^2 +  \al/4 < \al/4 + \al/4 = \al/2,
\label{recineq2}
\end{equation}
and
$$\exp(-3 n \phi(p)) < (\al/4)^3 < (\al/2)/4.$$

Hence, \eqref{conditionii} with $n$ replaced by $3 n$, and $\al$ replaced by $\al/2$ holds.
So we can iterate \eqref{recineq2}
and conclude that, for all integers $k \geq 0$, \\
$h(3,\,3^k n) < \alpha / (2^k)$.

The last part of the argument is exactly as for ordinary percolation: Note that if none of the reactangles
$[0, 3^{2 k +1} n] \times [0, 3^{2 k} n]$ and $[0, 3^{2 k +1} n] \times [0, 3^{2 k + 2} n]$,
$k=0, 1, 2,  \cdots$ has a white *-crossing in the `easy' (short) direction, then each of these rectangles has 
a black crossing in the long direction. Moreover, all these black crossings together form an infinite
occupied path.
Hence,
$$\th(p,\de) \geq 1 - \sum_{k = 0}^{\infty} h(3, n 3^k) \geq 1 - \alpha \sum_{k=0}^{\infty} (1/2)^k =  1- 2 \alpha >0.$$ 
This proves that (ii) implies (i).

\smallskip\noindent
Now we show that (i) implies (ii): Suppose $\th(p,\de) > 0$. Then there is (a.s.) an infinite occupied cluster, and by
Lemma \ref{unique} this cluster is unique. From the usual spatial symmetries, positive association, and the above mentioned
uniqueness one can, in exactly the same way as for ordinary percolation (see \cite{Gr}, Theorem 8.97) show that
$f(1,n) \ra 1$ as $n \ra \infty$. By Theorem \ref{thm-rsw-sdp} it follows that also $f(3,n) \ra 1$ 
as $n \ra \infty$; so (ii) holds.
\end{proof}
\label{fscriterion}
\end{section}
\begin{section}{Proof of Theorem \ref{mainthm}}
We are now ready to prove Theorem \ref{mainthm}:

\begin{proof}
For $p < p_c$, we have (see Section 1) $\th(p,\de) = \th(p+(1-p)\de)$, so that continuity follows from continuity for
ordinary percolation. If $p = p_c$ and $\de > p_c + \varepsilon$ for some $\varepsilon > 0$, then (trivially)
there is a neighborhood of $(p,\de)$ where the
sdp model dominates ordinary percolation with parameter $p_c+\varepsilon/2 > p_c$; hence $\th(\cdot,\cdot) > 0$ on this
neighborhood, and Proposition \ref{contprop} implies continuity of $\th(\cdot,\cdot)$ at $(p_c,\de)$.
(In fact, by combining this argument with an Aizenman-Grimmett type argument, one can extend this result and show
that there is an $\varepsilon >0$ such that $\th(\cdot,\cdot)$ is continuous at $(p_c,\de)$ if
$\de > p_c - \varepsilon$).

Finally, we consider the case where $p > p_c$. If $\th(p,\de) = 0$, continuity at $(p,\de)$ follows from
part (b) of Proposition \ref{contprop}.
So suppose $\th(p,\de) > 0$.
Let $\al$ as in Theorem \ref{fsizethm}.
By that theorem there is an $n \geq \hN(p)$ with
$$f_{p,\de}(3, n) > 1 -\al.$$
Hence, by Lemma \ref{cylcont} there is an open neighborhood $W$ of $(p,\de)$ such that
\begin{equation}
f_{p',\de'}(3, n) > 1 -\al,
\label{last}
\end{equation}
for all $(p',\de') \in W$.
Since $n \geq \hN(p)$ and $\hN(.)$ is decreasing, it follows from \eqref{last} and Theorem \ref{fsizethm}
that 
$\th(\cdot,\cdot) > 0$ on $S$,
%\begin{equation}
%\forall (p',\de') \in S \, \,\exists m \geq \hN(p') \mbox{ such that } f_{p',\de'}(3, m) > 1 -\al, 
%\label{gebied1}
%\end{equation}
where $S$ is the set of all 
$(p',\de') \in W$ with $p' \geq p$. 
From this and Lemma \ref{mon} we conclude that $\th(\cdot,\cdot)$ is also strictly positive
on the set 
$$U := \{(p',\de') \,: \, p' < p \mbox{ and } p'+(1-p')\de' > r+(1-r) \be \mbox{ for some } (r,\be) \in S\}.$$
It is easy to see that $S \cup U$ contains an open neighborhood of $(p,\de)$. Now it follows from
part (a) of Proposition \ref{contprop} that $\th(\cdot)$ is continuous at $(p,\de)$.
This completes the proof of the main theorem.
\end{proof}

\begin{rem}
A crucial role in the proof is the finite-size criterion, Theorem \ref{fsizethm}. That theorem has been
formulated for $p >p_c$.  When $p=p_c$ (or $< p_c$) the sdp model is an ordinary percolation model, for which
a similar criterion is known. In fact, for ordinary percolation we do not have the dependency problems which
led to the introduction of $\hN$. Consequently, for $p=p_c$ we can take $\hN = 1$. But, on the other hand,
if we let $p$ tend to $p_c$ from above, the upper bound on $\hN(p)$ obtained from our arguments in Section
\ref{fscriterion} tends to $\infty$. And that, in turn, comes from the fact that our bound on dependencies,
Lemma \ref{asindep}, is in terms of path probabilities for an ordinary percolation model (on the matching lattice,
with parameter $1-p$) which is subcritical
but approaches criticality (which makes these bounds worse and worse) as $p$ approaches $p_c$ from above.
This is essentially why the proof of Theorem \ref{mainthm} does not work
at $p_c$. Of course, if it {\it would} work, the conjecture referred to in Section 1 would be false. We hope that attempts
to stretch the arguments in our paper as far as possible will substantially increase insight in the conjecture and help to
obtain a solution.
\label{remfinal}
\end{rem}

\label{sec6}
\end{section}

\end{document}